\documentclass [11pt]{amsart}
\usepackage {amsmath, amssymb, amscd, graphicx, color, url, epsf}
\usepackage[all,knot,arc]{xy}

\newtheorem {theorem}{Theorem}[section]

\newtheorem {lemma}[theorem]{Lemma}
\newtheorem {claim}{Claim}
\newtheorem {proposition}[theorem]{Proposition}
\newtheorem {corollary}[theorem]{Corollary}

\newcommand{\abs}[1]{\lvert#1\rvert}

\newcommand\R{\mathbb{R}}

\newcommand\Ta{\mathbb{T}_\alpha}
\newcommand\Tb{\mathbb{T}_\beta}
\newcommand\Tg{\mathbb{T}_\gamma}
\newcommand\Tbp{\mathbb{T}_{\beta^{\prime}}}

\def\Sym{\mathrm{Sym}}

\newcommand{\HFh}{\widehat{HF}}
\newcommand{\HFm}{{HF}^-}
\newcommand{\CFm}{{CF}^-}
\newcommand{\CFh}{\widehat{CF}}
\def\CFKm{{CFK}^-}

\def\HFKh{{\widehat{\mathrm{HFK}}}}
\def\HFKm{\mathrm{HFK}^-}

\def\lambdap{\lambda_+}
\def\lambdam{\lambda_-}
\def\lambdaps{\lambda^S_+}

\def\lambdaph{\widehat{\lambda}_+}
\def\thetah{\widehat{\theta}}
\def\lambdamh{\widehat{\lambda}_-}
\def\xp{\mathbf{x}_+}
\def\xps{\mathbf{x}_+^S}
\def\xm{\mathbf{x}_-}
\def\x{\mathbf{x}}
\def\y{\mathbf{y}}
\def\F{\mathbb{Z}_2}

\def\CFKm{{CFK}^-}
\def\CFKh{\widehat{CFK}}

\def\P{\mathcal{P}}

\def\x{\mathbf x}
\def\y{\mathbf y}

\newcommand\D{\mathcal D}
\newcommand\alphas{\mbox{\boldmath$\alpha$}}
\newcommand\betas{\mbox{\boldmath$\beta$}}
\newcommand\betasp{\mbox{\boldmath$\beta^{\prime}$}}
\newcommand\gammas{\mbox{\boldmath$\gamma$}}

\newcommand\ws{\mathbf w}
\newcommand\zs{\mathbf z}

\begin{document}

\title[Transversely Non Simple Knots]{Transversely Non Simple Knots}

\author [Vera V\'ertesi]{Vera V\'ertesi}
\thanks {The author was supported by NSF grant number FRG-0244663 and OTKA 67867 and 67870.}
\address {Institute of Mathematics , E\"otv\"os Lor\'and University, Budapest, Hungary}
\email {wera@szit.bme.hu}

\begin {abstract}
By proving a connected sum formula for the Legendrian invariant $\lambdap$ in knot Floer homology
we exhibit infinitely many transversely non simple knots.
\end {abstract}

\maketitle
\section {Introduction}

The study of Legendrian and transverse knots is central in contact geometry.
A Legendrian knot with a given knot type has two classical invariants: the Thurston-Bennequin number and the rotation number. The problem of classifying Legendrian knots up to Legendrian isotopy naturally leads to the question whether these invariants classify Legendrian knots. A knot type is called Legendrian simple if any two realization of it with equal classical invariants are Legendrian isotopic.
For transverse knots there is only one classical invariant the self-linking number. Similarly, those knot types that, are classified by the self-linking number are called transversely simple.
The unknot \cite{EF}, torus knots and the figure-eight knot \cite{EH} were proved to be both Legendrian and transversely simple. By constructing a new invariant for Legendrian knots, Chekanov \cite{Ch} showed that not all knots are Legendrian simple, in particular he proved that the knot $5_2$ is not Legendrian simple. Later many other Legendrian non simple knots were detected by Epstein, Fuchs and Meyer \cite{EFM}, and by Ng \cite{Ng}.
The case for transverse knots turned out to be harder, Birman and Menasco \cite{BM}, and Etnyre and Honda \cite{EH} constructed families of transversely non simple knots using braid and convex surface theory. Recently Ng, Ozsv\'ath and Thurston \cite{NgOT} gave such examples using the Legendrian invariant in knot Floer homology.

Heegaard Floer homology $\HFh(Y)$, $\HFm(Y)$ defined by Ozsv\'ath and Szab\'o \cite{OSZ3m} are invariants for three-manifolds. The construction was extended \cite{OSZknot} to give the invariants $\HFKh(Y,K)$, $\HFKm(Y,K)$ for null-homologous knots $K\subset Y$ via doubly pointed Heegaard diagrams.
Using Heegaard diagrams with multiple basepoints the invariants were generalized for links \cite{AAA}. Multiply pointed Heegaard diagrams  turned out to be extremely useful in the case of knots too, and led to the discovery of a combinatorial version of the knot Floer homologies through grid diagrams \cite{MOSzT,MOS}. This version provided a natural way to define invariants $\lambdap$ and $\lambdam$ of Legendrian and $\theta$ for transverse knots
in the three-sphere \cite{OSzT}.

In this paper we will prove

\begin{theorem}\label{thm:connectedsum} Let $L_1$ and $L_2$ be (oriented) Legendrian knots in the standard contact 3-sphere of topological type $K_1$ and $K_2$. Then
there is an isomorphism:
$\HFKm(S^3,m(K_1))\otimes_{\F[U]}\HFKm(S^3,m(K_2)) \to \HFKm(S^3,m(K_1\#K_2))$
which maps $\lambdap(L_1)\otimes\lambdap(L_2)$ to $\lambdap(L_1\#L_2)$.
Similar statement holds for the $\lambdam$-invariant.
\end{theorem}

\begin{corollary}\label{cor:connectedsumhat}
Let $L_1$ and $L_2$ be (oriented) Legendrian knots in the standard contact 3-sphere of topological type $K_1$ and $K_2$. Then
there is an isomorphism:
$\HFKh(S^3,m(K_1))\otimes_{\F}\HFKh(S^3,m(K_2)) \to \HFKh(S^3,m(K_1\#K_2))$
which maps $\lambdaph(L_1)\otimes\lambdaph(L_2)$ to $\lambdaph(L_1\#L_2)$.
Similar statement holds for the $\lambdamh$-invariant.
\end{corollary}

Similar results hold for the $\theta$-invariant of transverse knots:

\begin{corollary}\label{thm:connectedsumtheta} Let $T_1$ and $T_2$ be transverse knots in the standard contact 3-sphere of topological type $K_1$ and $K_2$. Then
there are isomorphisms:
$\HFKm(S^3,m(K_1))\otimes_{\F[U]}\HFKm(S^3,m(K_2)) \to \HFKm(S^3,m(K_1\#K_2))$ and
$\HFKh(S^3,m(K_1))\otimes_{\F}\HFKh(S^3,m(K_2)) \to \HFKh(S^3,m(K_1\#K_2))$
which maps $\theta(T_1)\otimes\theta(T_2)$ to $\theta(T_1\#T_2)$
and $\thetah(T_1)\otimes\thetah(T_2)$ to $\thetah(T_1\#T_2)$, respectively.
\end{corollary}

As an application of the above result we prove (cf.\ also \cite{BM}):

\begin{theorem}\label{thm:nontransversesimple}
There exist infinitely many transversely non simple knots.
\end{theorem}

The above result can also be obtained by the connected sum formula of Honda and Etnyre \cite{EH2}.

%
%

The paper is organized as follows. In Section \ref{sec:preliminaries} we recall the definitions and collect the basic facts about Legendrian and transverse knots, knot Floer homology and the Legendrian and transverse invariant.
In Section \ref{sec:proofmain} we introduce spherical grid diagrams and prove Theorem \ref{thm:connectedsum}.
In Section \ref{sec:prooftransversenonsimple} we use the results of Section \ref{sec:proofmain} to prove
Theorems \ref{thm:nontransversesimple}. 
\vspace{0.2cm}
\begin{center}
\textsc{Acknowledgment}
\end{center}
\vspace{0.2cm}
I would like to thank Peter Ozsv\'ath and Andr\'as Stipsicz for for their guidance and help during the course of this work.
In particular several key ideas in the proof were suggested by Andr\'as Stipsicz. This work was done while I visited Columbia University, I am grateful for their hospitality. I also wish to thank John Baldwin for pointing out a mistake in the paper, and to Tam\'as Terpai for reading the manuscript of the paper.

\section{Preliminaries}\label{sec:preliminaries}

\subsection{Legendrian and transverse knots}

A \emph{Legendrian knot} $L$ in $\R^3$ (or in $S^3=\R^3\cup\{\infty\}$), endowed with the standard contact form $dz-ydx$ is an oriented knot along which the form  $dz-ydx$ vanishes identically.
Legendrian knots are determined by their front projection to the $xz$-plane; the projections are smooth in all but finitely many cusp points, have no vertical tangents ($y=\frac{dz}{dx}$) and at each crossing the strand with smaller slope is in the front. Note, that in order to have the standard orientation on $\R^3$ the $y$ axis points into the page. By changing the parts with vertical tangents to cusps and adding zig-zags, a generic smooth projection of a knot can be arranged to be of the above type. Thus any knot can be placed in Legendrian position. But this can be done in many different ways up to Legendrian isotopy. For example by adding extra zig-zags in the front projection we obtain a different Legendrian representative. This method is called \emph{stabilization}. If we add a down cusp then it is called \emph{positive stabilization} and it is called \emph{negative stabilization} otherwise. (Here and throughout the paper we use the conventions of Etnyre \cite{E}.) There are two classical invariants for Legendrian knots, described in the following. By pushing off the knot in the $\frac{\partial}{\partial z}$ direction we obtain the \emph{Thurston-Bennequin framing} of the Legendrian knot. Comparing this with the Seifert framing we get a number which is called the \emph{Thurston-Bennequin number} $\textrm{tb}(L)$. The \emph{rotation number} $\textrm{r}(L)$ is the winding number of $TL$ with respect to a trivialization of the contact planes along $L$, that extends to a Seifert surface.

A \emph{transverse knot} in $S^3$ with the standard contact structure is a knot along which the contact form $dz-ydx$ never vanishes.
Any transverse knot is naturally endowed with an orientation, the one along which the contact form is positive.
Again, every knot can be placed in transverse position, by for example translating its Legendrian realization in the $\pm\frac{\partial}{\partial y}$ direction. The resulting transverse knot is called the \emph{transverse push off} of the Legendrian knot. A push off is called positive if the orientation of the knot agrees with the natural orientation of the transverse knot and called negative otherwise.
Note that the negative push off is a transverse knot with reverse orientation. A Legendrian knot is called the \emph{Legendrian approximation} of its positive push off.
Two transverse knots are transverse isotopic if and only if their Legendrian approximations have common negative stabilizations.
By pushing off the transverse knot $T$ in a direction of a vector field in the contact planes we get a
framing. Comparing this to the Seifert framing we get the \emph{self-linking number} $\textrm{sl}(T)$. The self-linking number can be deduced from the classical invariants of a Legendrian approximation of the knot: $\textrm{sl}(L_{\pm})=\textrm{tb}(L)\mp\textrm{r}(L)$.

A knot is called \emph{Legendrian simple} (or \emph{transversely simple}) if any two Legendrian (transverse) realization of it with equal Thurston-Bennequin and rotation (self-linking) number(s) are isotopic through Legendrian (transverse) knots.

\begin{figure}
\centering
\includegraphics[scale=0.4]{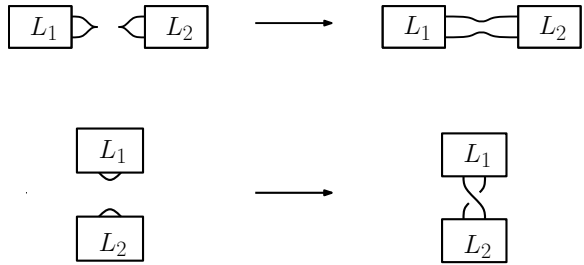}
\caption{Connected sum of two Legendrian knots}
\label{fig:legendreonnsum}
\end{figure}

As it is explained in \cite{EH2}, there is a well defined notion of the connected sum of two Legendrian or transverse knots in $S^3$, which comes from connected summing the two $S^3$'s the knots are sitting in. The above process has a good description in the front projection of the knots, as can be seen in Figure \ref{fig:legendreonnsum}.

\subsection{Knot Floer homology with multiple basepoints}

Here we outline the basic definitions of knot Floer homologies with multiple basepoints, originally defined by Ozsv\'ath and Szab\'o \cite{AAA} and independently by Rasmussen \cite{R}. Consider  a knot $K$ in an oriented, closed three-manifold $Y$, there is a self-indexing Morse function with $k$ critical points such that $K$
is made out of $2k$ flow lines connecting all the index zero and index three critical points. Such a Morse function gives rise to a Heegaard diagram $(\Sigma,\alphas,\betas,\ws,\zs)$ for $(Y,K)$ in the following way. Let $\Sigma=f^{-1}(\frac{3}{2})$ be a genus $g$ surface. The $\alpha$-curves $\alphas=\{\alpha_i\}_{i=1}^{g+k-1}$ are defined to be the circles of $\Sigma$ whose points flow down to the index one critical points. Similarly $\betas=\{\beta_i\}_{i=1}^{g+k-1}$ are the curves with points flowing up to the index two critical points. Finally let $\ws=\{w_i\}_{i=1}^{k}$ be the positive and $\zs=\{z_i\}_{i=1}^{k}$ the negative intersection points of $K$ with $\Sigma$.

Consider the module $\CFm(\Sigma, \alphas,\betas,\ws)$ over the polynomial algebra $\F[U_1,\dots,U_k]$ freely generated by the intersection points of the totally real submanifolds $\Ta=\alpha_1\times\cdots\times\alpha_{g+k-1}$ and $\Tb=\beta_1\times\cdots\times\beta_{g+k-1}$ of $\Sym^{g+k-1}(\Sigma)$. This module is endowed with the differential:
 \[
\partial^-\x=
\sum_{\mathbf{y}\in \Ta\cap\Tb}
\sum_{\begin{subarray}{l}
\phi \in\pi_2(\x,\mathbf{y}) \\
\mu(\phi)=1
\end{subarray}}
\left|\frac{\mathcal{M}(\phi)}{\R}\right|U_1^{n_{w_1}(\phi)}\cdots U_k^{n_{w_k}(\phi)}\mathbf{y}
\]
where as usual $\pi_2(\x,\y)$ denotes the space of homotopy classes of Whitney disks connecting $\x$ to $\mathbf{y}$, $\mathcal{M}(\phi)$ denotes the moduli space of pseudo-holomorphic representatives of $\phi$, the Maslov index,
$\mu(\phi)$ denotes its formal dimension and $n_p(\phi)=\#\{\phi^{-1}(p\times\Sym^{g+k-2}(\Sigma))\}$ is the local multiplicity of $\phi$ at $p$. Let

\begin{equation}\label{eqn:fact}
\left(\CFh(\Sigma,\alphas,\betas,\ws),\widehat{\partial}\right)=\left(\frac{\CFm(\Sigma,\alphas,\betas,\ws)}{(U_1=0)},\left[\partial^-\right]\right).
\end{equation}

The chain-homotopy type of the above complexes are invariants of $Y$ in the following sense:

\begin{theorem}\cite{AAA} Let $Y$ be an oriented closed three manifold. Consider the Heegaard diagrams $(\Sigma_1,\alphas_1,\betas_1,\ws_1)$ and $(\Sigma_2,\alphas_2,\betas_2,\ws_2)$ for $Y$ with $\abs{\ws_1}=k_1$ and $\abs{\ws_2}=k_2$. Assuming $k_1\geq k_2$ both complexes  $\CFm(\Sigma_1,\alphas_1,\betas_1,\ws_1)$ and $\CFm(\Sigma_2,\alphas_2,\betas_2,\ws_2)$ can be thought of as chain complexes over $\F[U_1,\dots,U_{k_1}]$ where $U_{k_2},\dots,U_{k_1}$ act identically on the latter complex. With this setup the two chain complexes are chain-homotopy equivalent.
Similar statement holds for the $\CFh$-theory, moreover the chain homotopy equivalences form a commutative diagram with the factorization map of (\ref{eqn:fact}). \qed
\end{theorem}

Hereafter we assume that our underlying three-manifold is the three-sphere. Note that in this case the homology of $\CFm(\Sigma,\alphas,\betas,\ws)$ is $\HFm(S^3)=\F[U]$.
The relative Maslov-grading of two intersection points $\x,\mathbf{y}\in \Ta\cap\Tb$ is defined by $M(\x)-M(\mathbf{y})=\mu(\phi)-2\sum n_{w_i}(\phi)$, where $\phi\in\pi_2(\x,\mathbf{y})$ is any homotopy class from $\x$ to $\mathbf{y}$. We extend this relative grading to the whole module by $M(U_1^{a_1}\cdots U_k^{a_k}\x)=M(\x)-2(a_1+\cdots+a_k)$. This grading can be lifted to an absolute grading in $S^3$ by fixing the grading of the generator of $\HFm(S^3)=\F[U]$ at $0$.

Note that so far we made no reference to the basepoints $\zs$. The relative Alexander grading is defined by $A(\x)-A(\mathbf{y})=\sum n_{z_i}(\phi)-\sum n_{w_i}(\phi)$, where again $\phi$ can be chosen to be any homotopy class in $\pi_2(\x,\mathbf{y})$. This relative grading can be uniquely lifted to an absolute Alexander grading to satisfy $\sum T^{A(x)}=\Delta_K(T)(1-T)^{n-1} \quad (\textrm{mod } 2)$, where $\Delta_K(T)$ is the symmetrized Alexander polynomial. We can extend the Alexander grading to the module by $A(U_1^{a_1}\cdots U_k^{a_k}\x)=A(\x)-(a_1+\cdots+a_k)$. As the local multiplicities of pseudo-holomorhic discs are non-negative, we obtain filtered chain complexes $\CFKm(\Sigma,\alphas,\betas,\ws,\zs)$ and $\CFKh(\Sigma,\alphas,\betas,\ws,\zs)$, that are invariants of the knot:

\begin{theorem}\cite{AAA} Let $K$ be an oriented knot in the 3-sphere. Consider the Heegaard diagrams $(\Sigma_1,\alphas_1,\betas_1,\ws_1,\zs_1)$ and $(\Sigma_2,\alphas_2,\betas_2,\ws_2,\zs_2)$ for $(S^3,K)$ with $\abs{\ws_1}=\abs{\zs_1}=k_1$ and $\abs{\ws_2}=\abs{\zs_2}=k_2$. Assuming $k_1\geq k_2$ both complexes $\CFKm(\Sigma_1,\alphas_1,\betas_1,\ws_1,\zs_1)$ and $\CFKm(\Sigma_2,\alphas_2,\betas_2,\ws_2,\zs_2)$ can be thought of as filtered chain complexes over $\F[U_1,\dots,U_{k_1}]$ where $U_{k_2},\dots,U_{k_1}$ act identically on the latter complex. With this setup the two filtered chain complexes are filtered chain-homotopy equivalent.
Similar statement holds for the $\CFKh$-theory, moreover the chain homotopy equivalences form a commutative diagram with the factorization map. \qed
\end{theorem}

As it is easier to work with we usually consider the homologies of the associated graded objects denoted by $\HFKm$. That is the homology of the complex $\CFKm(\Sigma,\alphas,\betas,\ws,\zs),\partial_0^-$ where,
\[
\partial^-_0\x=
\sum_{\mathbf{y}\in \Ta\cap\Tb}
\sum_{\begin{subarray}{l}
\phi \in\pi_2(\x,\mathbf{y}) \\
n_{z_1}(\phi)+\cdots+n_{z_k}(\phi)=0\\
\mu(\phi)=1
\end{subarray}}
\left|\frac{\mathcal{M}(\phi)}{\R}\right|U_1^{n_{w_1}(\phi)}\cdots U_k^{n_{w_k}(\phi)}\mathbf{y}.
\]

\begin{theorem}\cite{AAA} Let $K$ be an oriented knot in the 3-sphere $S^3$. Consider the Heegaard diagrams $(\Sigma_1,\alphas_1,\betas_1,\ws_1,\zs_1)$ and $(\Sigma_2,\alphas_2,\betas_2,\ws_2,\zs_2)$ for $(S^3,K)$. Then the $U_i$'s act identically on the homologies thus we can think of $\HFKm(\Sigma_1,\alphas_1,\betas_1,\ws_1,\zs_1)$ and $\HFKm(\Sigma_2,\alphas_2,\betas_2,\ws_2,\zs_2)$ as modules over $\F[U]$. And in this sense they are isomorphic.
Similar statement holds for the $\HFKh$-theory, moreover the isomorphisms form a commutative diagram with the factorization map. \qed
\end{theorem}


Knot Floer homology satisfies a K\"unneth principle for connected sums:

\begin{theorem}\cite{AAA,}\label{thm:connsum} Let $K_1$ and $K_2$ be oriented knots in $S^3$ described by the Heegaard diagrams $(\Sigma_1,\alphas_1,\betas_1,\ws_1,\zs_1)$ and $(\Sigma_2,\alphas_2,\betas_2,\ws_2,\zs_2)$. Let $w_1^1\in \ws_1$ and $z_1^2\in \zs_2$. Then
\begin{enumerate}
\item[(1)] $\left(\Sigma_1\#\Sigma_2,\alphas_1\cup\alphas_2,\betas_1\cup\betas_2,(\ws_1-w_1^1)\cup\ws_2,\zs_1\cup(\zs_2-z_1^2)\right)$ is a Heegaard diagram for $K_1\# K_2$. Here the connected sum $\Sigma_1\#\Sigma_2$ is taken in the regions containing $w_1^1\in\Sigma_1$ and $z_1^2\in \Sigma_2$;
\end{enumerate}
\noindent Let $\abs{\ws_1}=\abs{\zs_1}=k_1$ and $\abs{\ws_2}=\abs{\zs_2}=k_2$. Both complexes $\CFKm(\Sigma_1,\alphas_1,\betas_1,\ws_1,\zs_1)$ and $\CFKm(\Sigma_2,\alphas_2,\betas_2,\ws_2,\zs_2)$ are $\F[U_1,\dots,U_{k_1},V_1,\dots,V_{k_2}]$-modules with $U_1,\dots,U_{k_1}$ acting trivially on the latter and $V_1,\dots,V_{k_2}$ acting trivially on the former. Then
\begin{enumerate}

\item[(2)] $\CFKm\left(\Sigma_1,\alphas_1,\betas_1,\ws_1,\zs_1\right)\otimes_{U_1=V_1}\CFKm\left(\Sigma_2,\alphas_2,\betas_2,\ws_2,\zs_2\right)$ is filtered chain homotopy equivalent to  $\CFKm\left(\Sigma_1\#\Sigma_2,\alphas_1\cup\alphas_2,\betas_1\cup\betas_2,(\ws_1-w_1^1)\cup\ws_2,\zs_1\cup(\zs_2-z_1^2)\right)$;


\item[(3)] $\HFKm(S^3,K_1\# K_2)$ is isomorphic to $\HFKm(S^3,K_1)\otimes\HFKm(S^3,K_2)$ and this isomorphism can be given by $\x_1 \otimes \x_2\mapsto( \x_1,\x_2) $ on the generators.

%
%

\end{enumerate}
Here, all tensor products are taken over $\F[U_1,\dots,U_{k_1},V_1,\dots,V_{k_2}]$.

Similar statement holds for the $\CFKh$-theory, moreover the chain homotopy equivalences form a commutative diagram with the factorization map. \qed

\end{theorem}

\subsection{Grid diagrams}

\begin{figure}
\centering
\includegraphics[scale=0.5]{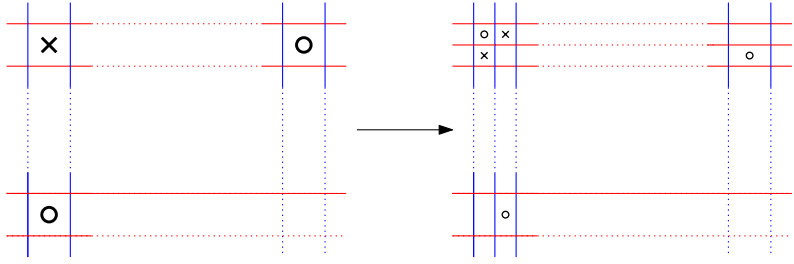}
\caption{Stabilization of type $X\!:\!NW$}
\label{fig:stab}
\end{figure}

As it was observed \cite{MOSzT, MOS} knot Floer homology admits a completely combinatorial description via grid diagrams.
A \emph{grid diagram} $G$ is a $k$ by $k$ square grid placed on the plane with some of its cells decorated with an $X$ or an $O$ and containing exactly one $X$ and $O$ in each of its rows and columns. Such a diagram naturally defines an oriented link projection by connecting the $O$'s to the $X$'s in each row and the $X$'s to the $O$'s in the columns and letting the vertical line to overpass at the intersection points. For simplicity we will assume that the corresponding link is a knot $K$. There are certain moves of the grid diagram that do not change the (topological) knot type \cite{OSzT}. These are
\emph{cyclic permutation} of the rows or columns, \emph{commutation} of two consecutive rows (columns) such that the $X$ and the $O$ from one row (column) does not separate the $X$ and the $O$ from the other row (column) and \emph{(de)stabilization} which is described in the following. A square in the grid containing an $X$ ($O$) can be subdivided into four squares by introducing a new vertical and a new horizontal line dividing the row and the column that contains that square. By replacing the $X$ ($O$) by one $O$ ($X$) and two $X$'s ($O$'s) in the diagonal of the new four squares and placing the two $O$'s ($X$'s) in the subdivided row and column appropriately we get a new grid diagram which is called the stabilization of the original one. The inverse of stabilization is destabilization. There are eight types of (de)stabilization: $O\!:\!SW$, $O\!:\!SE$, $O\!:\!NW$, $O\!:\!NE$, $X\!:\!SW$, $X\!:\!SE$, $X\!:\!NW$ and $X\!:\!NE$, where the first coordinate indicates which symbol we started with and the second shows the placement of the unique new symbol. A stabilization of type $X\!:\!NW$ is depicted on Figure \ref{fig:stab}.

Placing the grid on a torus by identifying the opposite edges of the square grid we obtain a Heegaard diagram with multiple basepoints for $(S^3,K)$. Here the $\ws$'s correspond to the $O$'s and the $\zs$'s to the $X$'s. As each region of this Heegaard diagram is a square, it is``nice'' in the sense defined in \cite{SW}. Thus boundary maps can be given by rectangles. This observation led to a completely combinatorial description of knot Floer homology \cite{MOSzT,MOS} in the three-sphere. This can be set up without referring to the original holomorphic theory \cite{MOS}.

\subsection{Legendrian and transverse invariants on grid diagrams}\label{subsec:leginv}

Consider a grid diagram $G$ it does not only describes a knot projection but also a front projection of a Legendrian realization of its mirror $m(K)$ as follows. Rotate the grid diagram by $45^{\circ}\!$ clockwise, reverse the over- and under crossings and turn the corners into cusps or smooth them as appropriate. Legendrian Reidemeister moves correspond to certain grid moves giving the following result:

\begin{proposition}\cite{OSzT}\label{prop:leggrid} Two grid diagrams represent the same Legendrian knot if and only if they can be connected by a sequence of cyclic permutation, commutation, and (de)stabilization of types $X\!:\!NW$, $X\!:\!SE$, $O\!:\!NW$ and $O\!:\!SE$. \qed
\end{proposition}

Moreover stabilization of type $X\!:\!NE$ or $O\!:\!SW$ of the grid diagram correspond to negative stabilization of the knot thus

\begin{proposition}\cite{OSzT}
Two grid diagram represent the same transverse knot if and only if they can be connected by a sequence of cyclic permutation, commutation, and (de)stabilization of types $X\!:\!NW$, $X\!:\!SE$, $X\!:\!NE$, $O\!:\!NW$, $O\!:\!SE$ and $O\!:\!SW$. \qed
\end{proposition}

Consider a grid diagram $G$ for a Legendrian knot $L$ of knot type $K$,
pick the upper right corner of every cell containing an $X$. This gives a generator of $\CFKm(S^3,m(K))$ denoted by $\xp(G)$. Here $m(K)$ denotes the mirror of $K$. Since there is no positive rectangle starting at $\xp(G)$, it is
a cycle defining an element $\lambdap(G)$ in $\HFKm(S^3,m(K))$. Similarly one can define $\xm(G)$ and $\lambdam(G)$ by taking the lower left corners of the cells containing $X$'s. These elements are proved to be independent of the grid moves that preserve the Legendrian knot type, giving an invariant of the Legendrian knot $L$:
\begin{theorem}\cite{OSzT} Consider two grid diagram $G_1$ and $G_2$ defining the same oriented Legendrian knot, then there is a quasi-isomorphism of the graded chain complexes $\CFKm$ taking $\xp(G_1)$ to $\xp(G_2)$ and $\xm(G_1)$ to $\xm(G_2)$. \qed
\end{theorem}

One can understand their transformation under positive and negative stabilization:

\begin{theorem}\cite{OSzT} Let $L$ be an oriented Legendrian knot, denote $L_+$ its positive and $L_-$ its negative stabilization. Then
\begin{enumerate}
\item There is a quasi-isomorphism of the corresponding graded complexes taking $\lambdap(L)$ to $\lambdap(L_+)$ and $U\lambdam(L)$ to $\lambdam(L_+)$;
\item There is a quasi-isomorphism of the corresponding graded complexes taking $U\lambdap(L)$ to $\lambdap(L_-)$ and $\lambdam(L)$ to $\lambdam(L_-)$.
\end{enumerate}  \qed
\end{theorem}

We know \cite{E} that the Legendrian knots with transversely isotopic positive push offs are the ones that obtained by a sequence of negative (de)stabilization of $L$.
So there is a well defined invariant for transverse knots: if $L$ is a Legendrian approximation of $T$ then define $\theta(T)=\lambdap(L)$.

\begin{theorem}\cite{OSzT}
For any two grid diagrams $G_1$ and $G_2$ of Legendrian approximations of the transverse knot $T$ there is  quasi-isomorphism of the corresponding graded chain complexes taking $\theta(G_1)$ to $\theta(G_2)$. \qed
\end{theorem}

\section{Proof of Theorem \ref{thm:connectedsum}}\label{sec:proofmain}

The Legendrian invariant has two appearance depending on in which version of the Floer homology we think it is.
The one introduced in subsection \ref{subsec:leginv} is in the combinatorial Heegaard Floer homology. Once the grid is placed on the torus we get a Heegaard diagram, and thus there is a natural identification of the combinatorial Heegaard Floer complex with the holomorphic Heegaard Floer complex, and under this identification the previously defined invariant has a counterpart in the original, holomorphic Heegaard Floer homology. We will use the same notation for both. We introduce yet another invariant for Legendrian knots:


\subsection{Legendrian invariant on sphere Heegaard diagrams}

\begin{figure}
\centering
$\begin{array}{c@{\hspace{2cm}}c}
\includegraphics[scale=0.2]{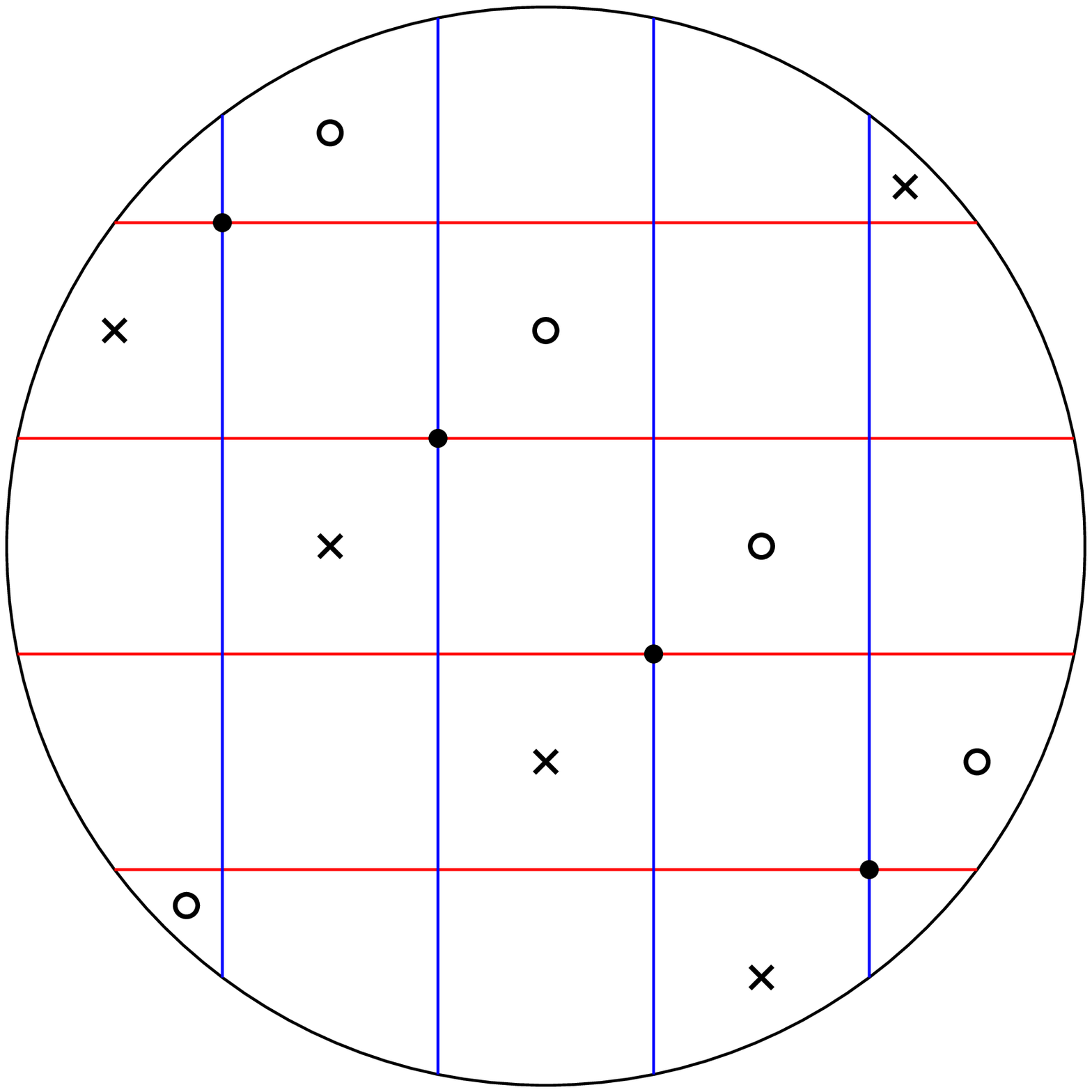}
&
\includegraphics[scale=0.2]{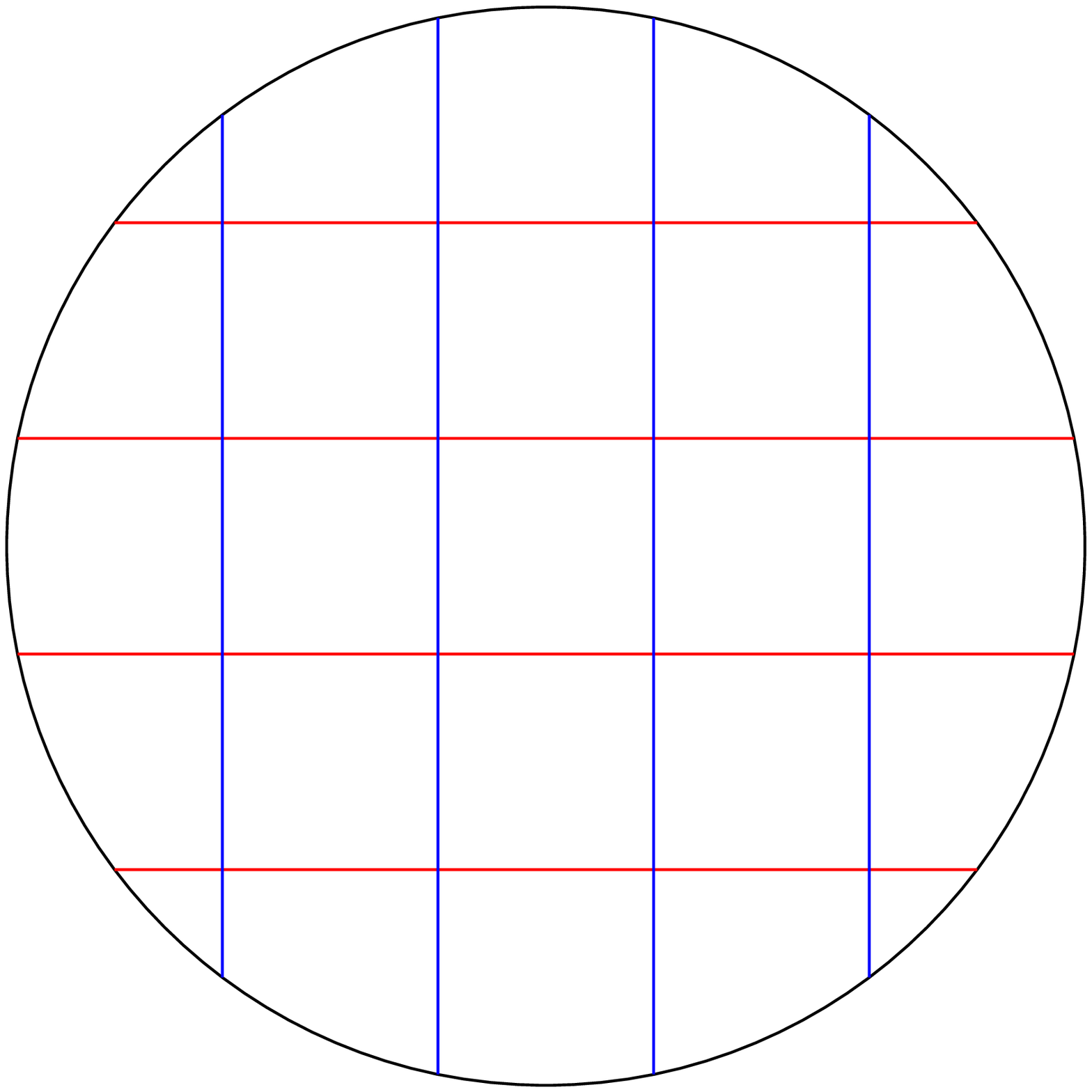}\\[0.2cm]
\mbox{front hemisphere} & \mbox{rear hemisphere}
\end{array}$
\caption{Spherical grid diagram for the trefoil knot}
\label{fig:trefoil}
\end{figure}

A $k$ by $k$ grid diagram $G$ of a Legendrian knot $L$ of topological type $K$ can also be placed on the two sphere in the following way.
Let $S^2=\{(x,y,z)\in\R^3:\; \vert(x,y,z)\vert=1 \}$ and define the circles $\widetilde{\alphas}=\{\widetilde{\alpha}_i\}_{i=1}^{k-1}$ as the intersection of $S^2$ with the planes $A_i=\{(x,y,z)\in\R^3:\; z=\frac{i}{k}-\frac{1}{2}\}$ ($i=1,\dots,k-1$), similarly $\widetilde{\betas}=\{\widetilde{\beta}_i\}_{i=1}^{k-1}$ as the intersection of $S^2$ with the planes $B_i=\{(x,y,z)\in\R^3:\; x=\frac{i}{k}-\frac{1}{2}
\}$ ($i=1,\dots,k-1$). Call $F=\{(x,y,z)\in\R^3:\; \vert(x,y,z)\vert=1, y\ge 0 \}$ the front hemisphere, and $R=\{(x,y,z)\in\R^3:\; \vert(x,y,z)\vert=1, y\le 0 \}$ the rear hemisphere. Then there is a grid on both the front and on the rear hemisphere. We place the $X$'s and the $O$'s on the front hemisphere in the way they were placed on the original grid $G$.
After identifying the $O$'s with $\widetilde{\ws}=\{\widetilde{w}_i\}_{i=1}^{k}$ and the $X$'s with $\widetilde{\zs}=\{\widetilde{z}_i\}_{i=1}^{k}$ this defines a Heegaard diagram with multiple basepoints $(S^2,\widetilde{\alphas},\widetilde{\betas},\widetilde{\ws},\widetilde{\zs})$ for $(S^3,K)$.
A spherical grid diagram for the trefoil knot can be seen on Figure \ref{fig:trefoil}.

Let $L$ be a Legendrian knot in $S^3$. To define the ``spherical'' Legendrian invariant $\lambdaps (L)$ we will use grid diagrams that have an $X$ in its upper right corner. This can always arranged by cyclic permutation, but in the following we will need a slightly stronger property:

\begin{lemma}\label{lem:ox} For any Legendrian knot there exists a grid diagram representing it which contains an $X$ in its upper right corner and an $O$ in its lower left corner.
\end{lemma}

\begin{proof}
Consider any grid diagram describing the Legendrian knot $L$.
As it is illustrated on Figure \ref{fig:ox} we can obtain a suitable diagram as follows.
First do a stabilization of type $X\!:\!NE$, and then do a stabilization of type $O\!:\!NE$ on the newly obtained $O$. Lastly by cyclic permutation we can place the lower $X$ introduced in the first stabilization to the upper right corner of the diagram. Notice that the $O$ on the upper right of this $X$ will be automatically placed to the lower left corner. According to Proposition \ref{prop:leggrid} the Legendrian type of the knot is fixed under these moves. Thus the statement follows.

\begin{figure}[h]
\centering
\includegraphics[scale=0.6]{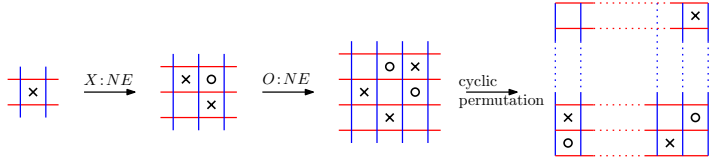}
\caption{Grid moves}
\label{fig:ox}
\end{figure}

\end{proof}

Suppose, that $G$ is a grid diagram having an $X$ in its upper right corner. Form a spherical grid diagram as above,
then $\xps(L)$ is the generator of $\CFKm(S^2,\widetilde{\alphas},\widetilde{\betas},\widetilde{\ws},\widetilde{\zs})$ consisting of those intersection points on the front hemisphere that occupy the upper right corner of each region marked with an $X$. Note that the $X$ in the upper right corner has no such corner. On Figure \ref{fig:trefoil} the element $\xps$ is indicated for the trefoil knot.  Similarly to the toroidal case:

\begin{lemma}\label{lem:xps}The element $\xps(L)$ is a cycle in $(S^2,\widetilde{\alphas},\widetilde{\betas},\widetilde{\ws},\widetilde{\zs})$.
\end{lemma}

\begin{proof}
We will show, that for any $\y$ there is no positive disc $\psi \in \pi_2(\xps,\y)$  with $\mu(\psi)=1$. As the diagram $\CFKm(S^2,\widetilde{\alphas},\widetilde{\betas},\widetilde{\ws},\widetilde{\zs})$ is ``nice'' in the sense of \cite{SW} the elements $\xps$ and $\y$ differ in exactly two coordinates and $\mathcal{D}(\psi)$ is either a rectangle or the union of two bigons. In any case $\mathcal{D}(\psi)$ contains an $X$ which means it is not counted in the boundary map.
\end{proof}

The homology class of $\xps$, denoted by $\lambdaps(G)$, turns out to be an invariant of $L$ (i.e.\ it is independent of the choice of the grid diagram, and the way it is placed on the sphere). This can be proved directly through grid moves but instead we show:

\begin{theorem}\label{thm:sphericalinv} Consider a grid diagram for  the Legendrian knot $L$ in $S^3$ having an $X$ in its upper right corner. Then there is a filtered quasi-isomorphism $\psi: \CFKm(T^2,\alphas,\betas,\ws,\zs) \to \CFKm(S^2,\widetilde{\alphas},\widetilde{\betas},\widetilde{\ws},\widetilde{\zs})$ of the corresponding toroidal and spherical Heegaard diagrams which maps $\xp(L)$ to $\xps(L)$.
\end{theorem}

In the proof we will need the notion of Heegaard triples, which we will briefly describe for completeness. Consider a pointed Heegaard triple $(\Sigma,\alphas,\betas,\gammas,\zs)$ then the pairs $(\Sigma,\alphas,\betas,\zs)$, $(\Sigma,\betas,\gammas,\zs)$ and $(\Sigma,\alphas,\gammas,\zs)$ define the
three-manifolds $Y_{\alpha\beta}$, $Y_{\gamma\beta}$ and $Y_{\alpha\gamma}$, respectively. There is a map from $\CFm(\Sigma,\alphas,\betas,\zs)\otimes\CFm(\Sigma,\betas,\gammas,\zs)$ to $\CFm(\Sigma,\alphas,\gammas,\zs)$ given on
a generator $\x\otimes\y$ by
 \[
\sum_{\mathbf{v}\in \Ta\cap\Tg}
\sum_{\begin{subarray}{l}
u\in\pi_2(\x,\y,\mathbf{v}) \\
n_\zs(u)=0\\
\mu(u)=0
\end{subarray}}
\abs{\mathcal{M}(u)}\mathbf{v}
\]
where $\pi_2(\x,\y,\mathbf{v})$ is the set of homotopy classes of holomorphic triangles connecting $\x$, $\y$ to $\mathbf{v}$; maps from a triangle to $\Sym^{g+k-1}(\Sigma)$ sending the edges of the triangle to $\Ta$,$\Tb$ and $\Tg$. This gives a well-defined map on the homologies $\HFm$. Also, the same definition gives a map on the filtered chain complexes $\CFKm$. When $\gammas$ can be obtained from $\betas$ by Heegaard moves then the manifold $Y_{\beta\gamma}$ is $\#^{g}S^1\times S^2$ and $\HFm(\#^{g}S^1\times S^2)$ is a free $\F[U]$-module generated by $2^g$-elements. Denote its top-generator by $\Theta^-_{\beta\gamma}$. The map $\CFKm(Y_{\alpha\beta})\to \CFKm(Y_{\alpha\gamma})$ sending $\x$ to the image of $\x\otimes\Theta^-_{\beta\gamma}$ defines a quasi-isomorphism of the chain complexes.

\begin{proof}[Proof of theorem \ref{thm:sphericalinv}]
{From} a toroidal grid diagram one can obtain a spherical one by first sliding every $\beta$-curve over $\beta_1$ and every $\alpha$-curve
over $\alpha_1$ and then destabilize the diagram at $\alpha_1$ and $\beta_1$. Thus we will construct the quasi-isomorphism by the composition $\psi=\psi_{\textrm{destab}}\circ\psi_\alpha\circ\psi_\beta$, where
 \[
\psi_{\beta}=
\sum_{\mathbf{y}\in \Ta\cap\Tbp}
\sum_{\begin{subarray}{l}
u\in\pi_2(\xp(L),\Theta^-,\mathbf{y}) \\
n_\zs(u)=0\\
\mu(u)=0
\end{subarray}}
\abs{\mathcal{M}(u)}\mathbf{y}\]
with $\Theta^-\in\Tb\cap\Tbp$ being the top generator of $\HFm(T^2,\betas,\betasp,\zs)=\HFm(S^1\times S^2)$, and $\psi_\alpha$ is defined similarly.
Note that in the case of the sliding there is also a ``closest point'' map denoted by $^{\prime}$ for the sliding of the $\beta$-curves and by $^{\prime\prime}$ for the sliding of the $\alpha$-curves. We claim:

\begin{lemma}\label{lem:beta} $\psi_{\beta}(\xp)=\xp^{\prime}$.
\end{lemma}

\begin{lemma}\label{lem:alpha} $\psi_{\alpha}(\xp^{\prime})=(\xp^{\prime})^{\prime\prime}$.
\end{lemma}

Here we just include the proof of Lemma \ref{lem:beta}; Lemma \ref{lem:alpha} follows similarly.

\begin{proof}[Proof of Lemma \ref{lem:beta}]
Figure \ref{fig:slidebeta} shows a weekly admissible diagram for the slides of the $\beta$-curves.
\begin{figure}
\centering
\includegraphics[scale=0.4]{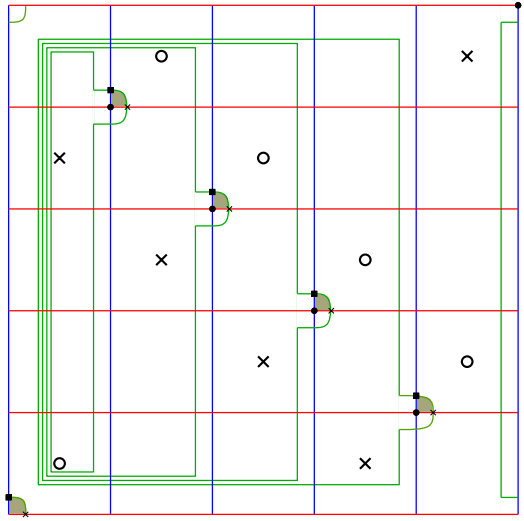}
\caption{Handleslides}
\label{fig:slidebeta}
\end{figure}

\begin{claim}
The Heegaard triple $(T^2,\alphas,\betas,\betasp,\zs)$ of Figure \ref{fig:slidebeta} is weekly admissible.
\end{claim}
\begin{proof}
Let $\mathcal{P}_{\beta_i \beta^{\prime}_i\beta_1}$ ($i>1$) denote the domain bounded by $\beta_i$, $\beta^{\prime}_i$ and $\beta_1$ and containing no basepoint, similarly $\mathcal{P}_{\beta_1\beta^{\prime}_1}$ denotes the domain bounded by $\beta_1$ and $\beta^{\prime}_1$ and containing no basepoint.
These domains form a basis for the periodic domains of $(T^2,\betas,\betasp,\zs)$ and as all have domains with both positive and negative coefficients we can see that $(T^2,\betas,\betasp,\zs)$ is weekly admissible.
Consider a triply periodic domain $\mathcal{P}$. If there is no $\alpha$-curve in its boundary, then it is a periodic domain of $(T^2,\betas,\betasp,\zs)$, and by the previous observation we are done. So $\mathcal{P}$ must contain an $\alpha$-curve in its boundary. To ensure it does not contain an $X$, there must be some vertical curve, either from $\beta$ or $\beta^{\prime}$, in the boundary. At the intersection point of the horizontal an vertical lines the domain must change sign.
\end{proof}

The grey area in Figure \ref{fig:slidebeta} indicates a domain of a canonical triangle $u_0$ in $\pi_2(\xp(L),\Theta^-,$ $\xp^{\prime}(L))$; by the
Riemann mapping theorem there is exactly one map with that domain. We claim that this is the only map that is encountered in $\psi_{\beta}$.
For this let $u\in \pi_2(\xp(L),\Theta^-,\mathbf{y})$ be a holomorphic triangle with  $\mu(u)=0$ and $n_{\zs}(u)=0$.
\begin{claim} There exists a periodic domain $\P_{\beta\beta^{\prime}}$ of $(S^2,\betas,\betasp,\zs)$ such that $\partial(\D(u)-\D(u_0)-\P_{\beta\beta^{\prime}})\vert_{\betas}=\emptyset$. Thus $(\D(u)-\D(u_0)-\P_{\beta\beta^{\prime}})\vert_{\betas}$ is a domain in $(T^2,\alphas,\betasp,\zs)$, representing an element $v$ in $\pi_2(\xp^{\prime},\mathbf{y})$ with Maslov index $\mu(v)=\mu(u)-\mu(u_0)-\mu(\P_{\beta\beta^{\prime}})=0$.
\end{claim}
\begin{proof}
As $n_{\zs}(u)=0$ and
$\xps(L)$ is in the upper right corner of the $X$'s, the domain of any triangle must contain $\mathcal{D}(u_0)$. Consequently
$\partial \mathcal{D}(u)\vert_{\beta_i}$ is an arc containing the small part $\overline{\mathcal{D}(u_0)}\cap \beta_i$ and some copies of the whole $\beta_i$. By subtracting $\mathcal{D}(u_0)$ and sufficiently many copies of the periodic domains $\mathcal{P}_{\beta_i \beta^{\prime}_i\beta_1}$ 
we obtain a domain with no boundary component on $\beta_i$. Doing the same process for every $i>1$ and then by subtracting some $\mathcal{P}_{\beta_1\beta_1^{\prime}}$ we can eliminate every $\beta_i$ from the boundary.
\end{proof}

\begin{claim} There is no positive disc in $\pi_2(\xp^{\prime},\mathbf{y})$.
\end{claim}
\begin{proof} This follows similarly as Lemma \ref{lem:xps}.
\end{proof}
\begin{claim}
None of the regions of $(T^2,\alphas,\betasp,\zs)$ can be covered completely with the periodic domains of $(S^2,\betas,\betasp,\zs)$ and $\D(u_0)$.
\end{claim}
\begin{proof} The periodic domains are the linear combinations of $\{\mathcal{P}_{\beta_i, \beta^{\prime}_i,\beta_1}\}_{i=2}^{k}\cup\{\mathcal{P}_{\beta_1, \beta^{\prime}_1,}\}$, and those can not cover the domains of $(S^2,\betas,\betasp,\zs)$.
\end{proof}
Putting these together we have, that $\D(u)-\D(u_0)-\P_{\beta\beta^{\prime}}$ have a negative coefficient, which gives a negative coefficient in $\D(u)$ too, contradicting the fact that $u$ was holomorphic. This proves Lemma \ref{lem:beta}.
\end{proof}

Note that by assuming that there is an $X$ in the upper right corner of the grid diagram we assured that the intersection point $\xp$ contains $\alpha_1\cap\beta_1$, and that point remained there. Thus by stabilizing at $\alpha_1$ and $\beta_1$ we get Theorem \ref{thm:sphericalinv}.
\end{proof}

\begin{proof}[Proof of Theorem \ref{thm:connectedsum}]

Consider two Legendrian knots $L_1$ and $L_2$ of topological type $K_1$ and $K_2$.
Note that once we obtain the result for $\lambdaps$, we are done. Indeed, passing from the toroidal diagram to the spherical one, the invariants $\lambdap(L_1)$ and $\lambdap(L_2)$ are mapped to $\lambdaps(L_1)$ and $\lambdaps(L_2)$, respectively. Knowing that $\lambdaps(L_1)\otimes\lambdaps(L_2)$ is mapped to $\lambdaps(L_1\#L_2)$ and passing back to the toroidal diagram there is an isomorphism that maps this to $\lambdap(L_1\#L_2)$. So the composition of the aboves proves Theorem \ref{thm:connectedsum}.

 Consider the grid diagrams $G_1$ and $G_2$ corresponding to $L_1$ and $L_2$ admitting the conditions of Lemma \ref{lem:ox}. These define the spherical grid diagrams
  $(S^2,\alphas_1,\betas_1,\ws_1,\zs_1)$ and $(S^2,\alphas_2,\betas_2,\ws_2,\zs_2)$. Let $z\in\zs_1$, $w\in\ws_2$ be the basepoints corresponding to the $X$ in the upper right corner of the first diagram and the $O$ in the lower left corner of the second diagram. Form the connected sum of
$(S^2,\alphas_1,\betas_1,\ws_1,\zs_1)$ and $(S^2,\alphas_2,\betas_2,\ws_2,\zs_2)$ at the regions containing $z$ and $w$ to obtain a Heegaard diagram with multiple basepoints $(S^2,\alphas_1\cup\alphas_2,\betas_1\cup\betas_2,\ws_1\cup(\ws_2-\{w\}),(\zs_1-\{z\}\cup\zs_2))$ of $(S^3,L_1\#L_2)$.
By \ref{thm:connsum} the map 
\begin{eqnarray*}
\psi_{\textrm{connsum}}:&&\HFKm(S^2,\alphas_1,\betas_1,\ws_1,\zs_1) \otimes \HFKm(S^2,\alphas_2,\betas_2,\ws_2,\zs_2) \to \\&&\HFKm(S^2,\alphas_1\cup\alphas_2,\betas_1\cup\betas_2,\ws_1\cup(\ws_2-\{w\}),(\zs_1-\{z\})\cup\zs_2)
\end{eqnarray*}
defined on the generators as $\mathbf{x}_1\otimes\mathbf{x}_2\mapsto (\mathbf{x}_1,\mathbf{x}_2)$ is an isomorphism. Thus the image of $\lambdaps(L_1)\otimes\lambdaps(L_2)$ is
$(\lambdaps(L_1),\lambdaps(L_2))$.

\begin{figure}
\centering
$\begin{array}{c@{\hspace{2cm}}c}
\includegraphics[scale=0.2]{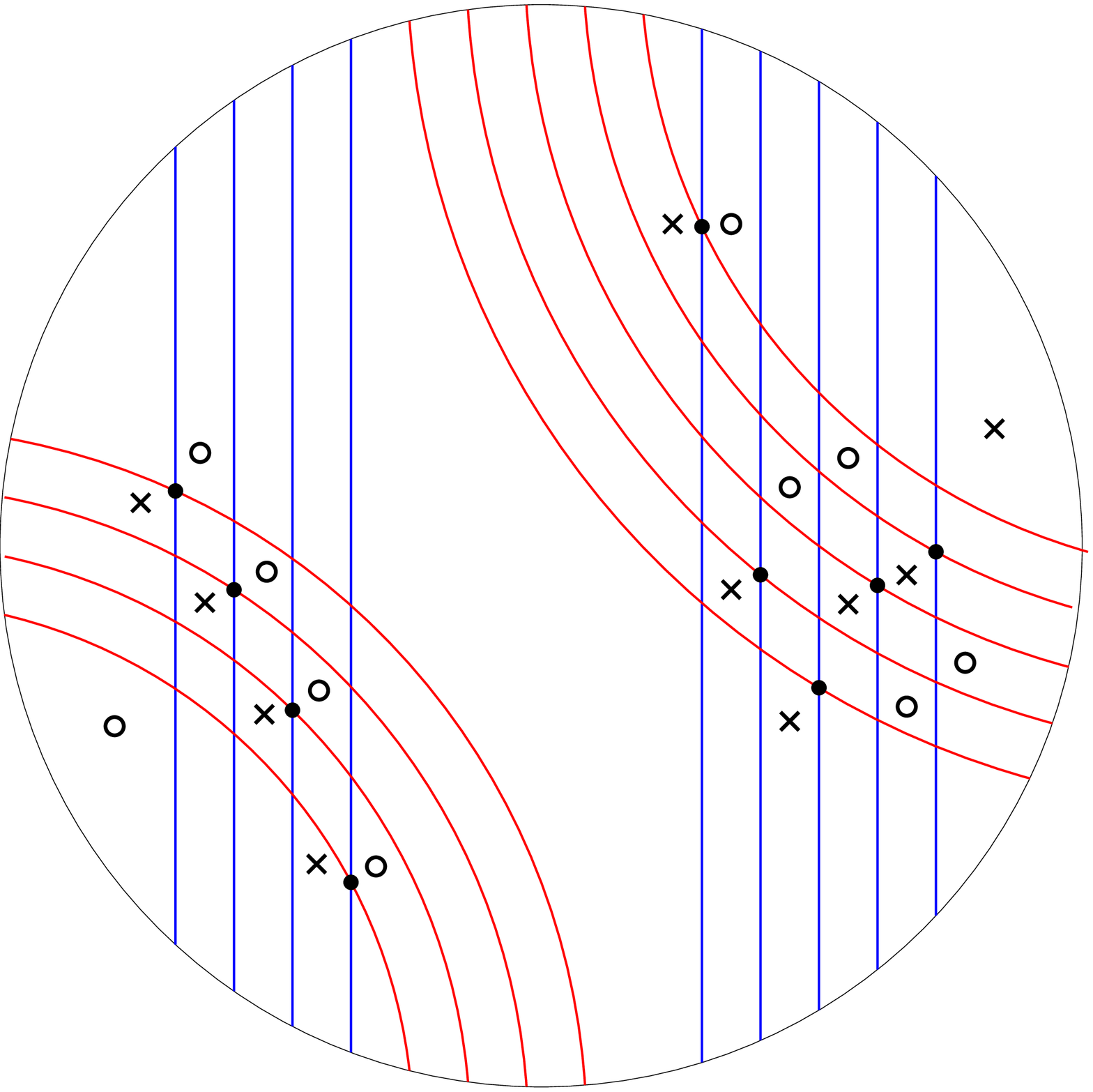}
&
\includegraphics[scale=0.2]{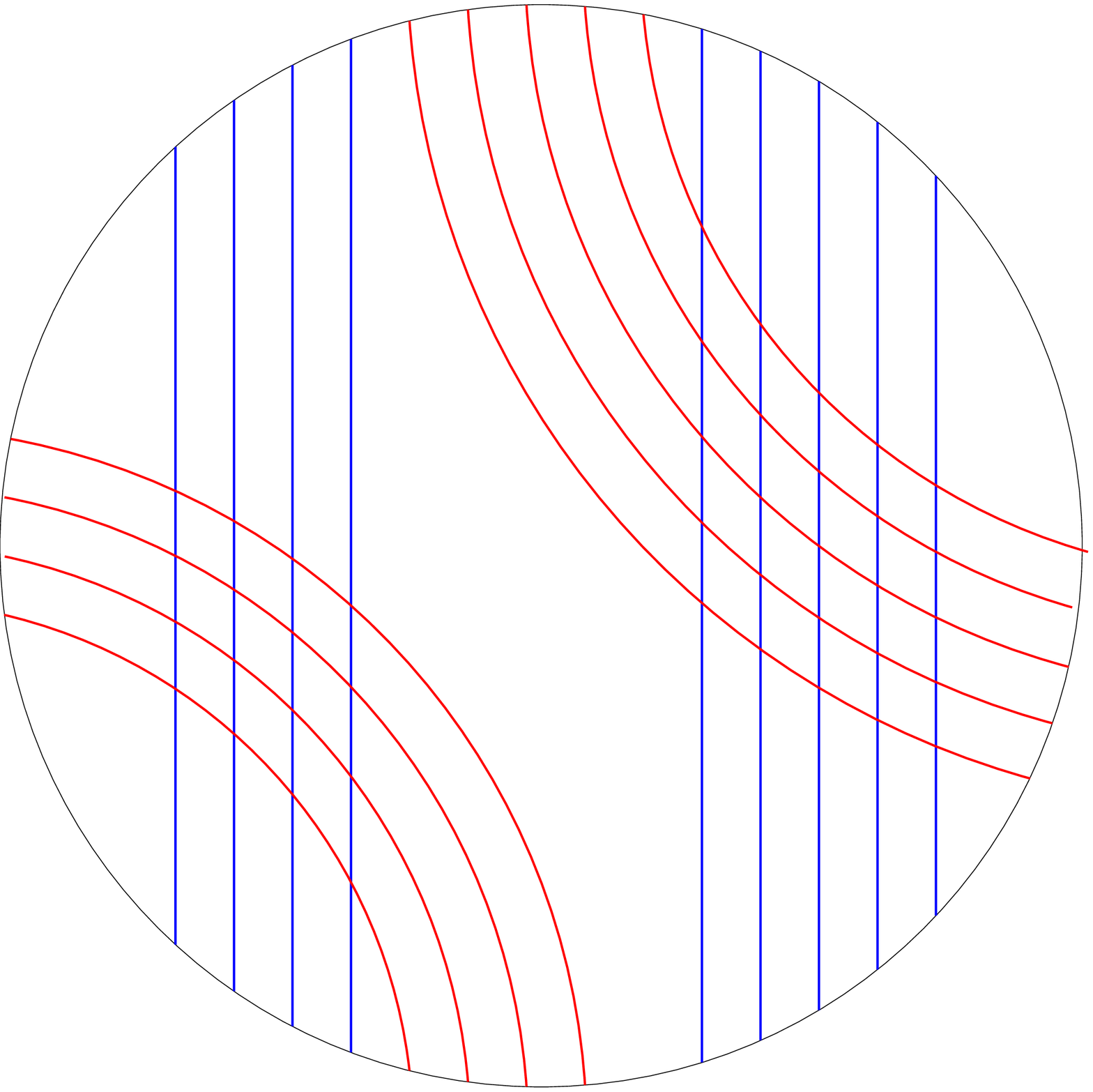}\\[0.2cm]
\mbox{front hemisphere} & \mbox{rear hemisphere}
\end{array}$
\caption{Connected Sum}
\label{fig:connectedsum}
\end{figure}

Figure \ref{fig:connectedsum} shows the resulting Heegaard diagram for the connected sum of a trefoil and a figure-eight knot.
From this diagram of the connected sum one can easily obtain a spherical grid diagram by isotoping every curve in $\alphas_1$ over the curves in $\betas_2$ and every curve in $\alphas_2$ over the curves in $\betas_1$ as shown on Figure \ref{fig:isotopy}. Indeed, the resulting diagram is a grid obtained by patching $G_1$ and $G_2$ together in the upper right $X$ of $G_1$ and the lower left $O$ of $G_2$ and deleting the $X$ and $O$ at issue. Now by connecting the $X$ in the lower row of $G_2$ to the $O$ in the upper row of $G_1$ and proceeding similarly in the columns, we get that the grid corresponds to the front projection of $L_1\#L_2$. Again, a quasi-isomorphism $\psi_{\textrm{isot}}$ is given  with the help of holomorphic triangles.
A similar argument as in the proof of Lemma \ref{lem:beta} shows that
under the isomorphism induced by $\psi_{\textrm{isot}}$ on the homologies the element $(\lambdaps(L_1),\lambdaps(L_2))$ is mapped to $\lambdaps(L_1\#L_2)$.

\begin{figure}
\centering
$\begin{array}{c@{\hspace{2cm}}c}
\includegraphics[scale=0.2]{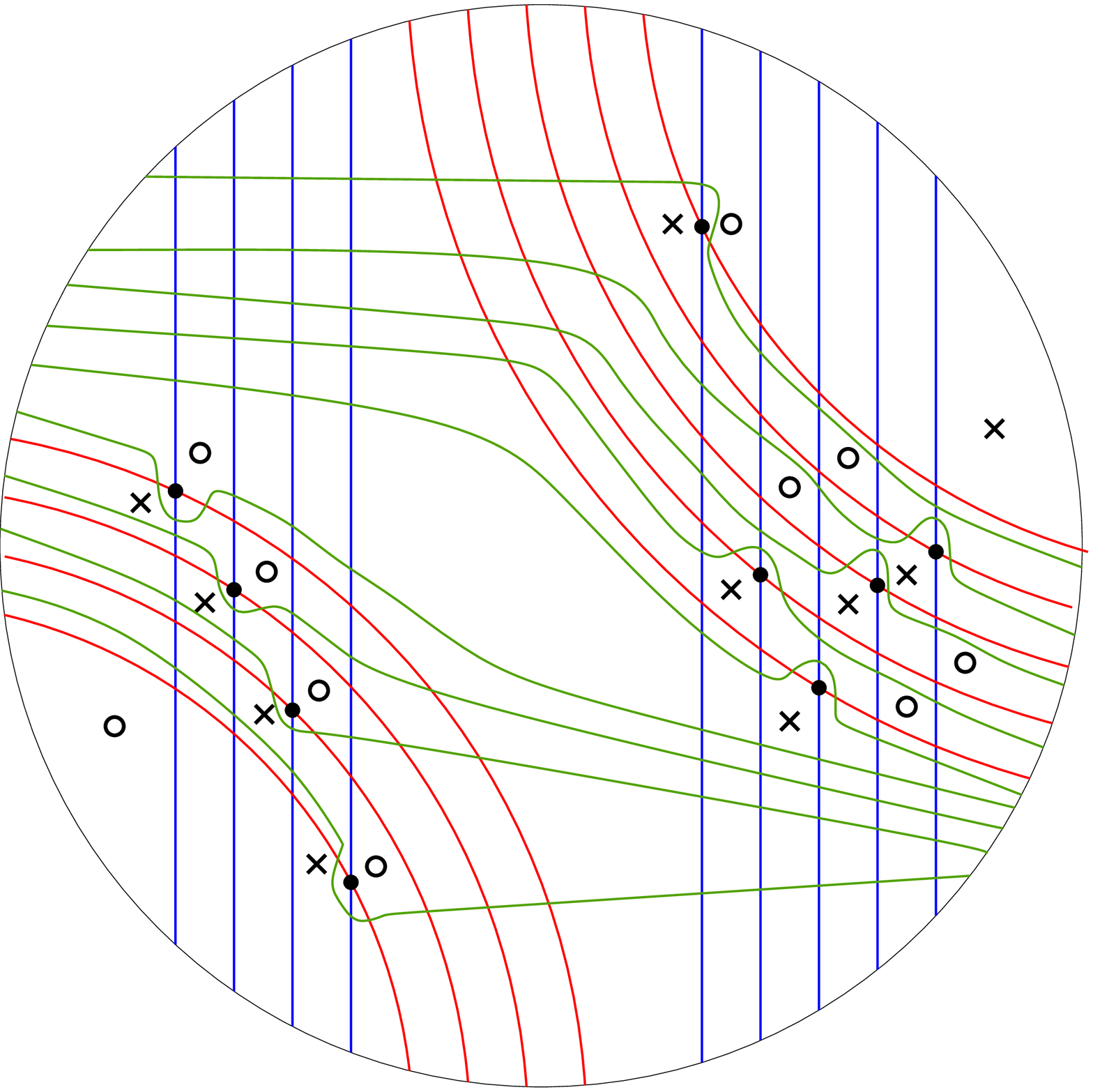}
&
\includegraphics[scale=0.2]{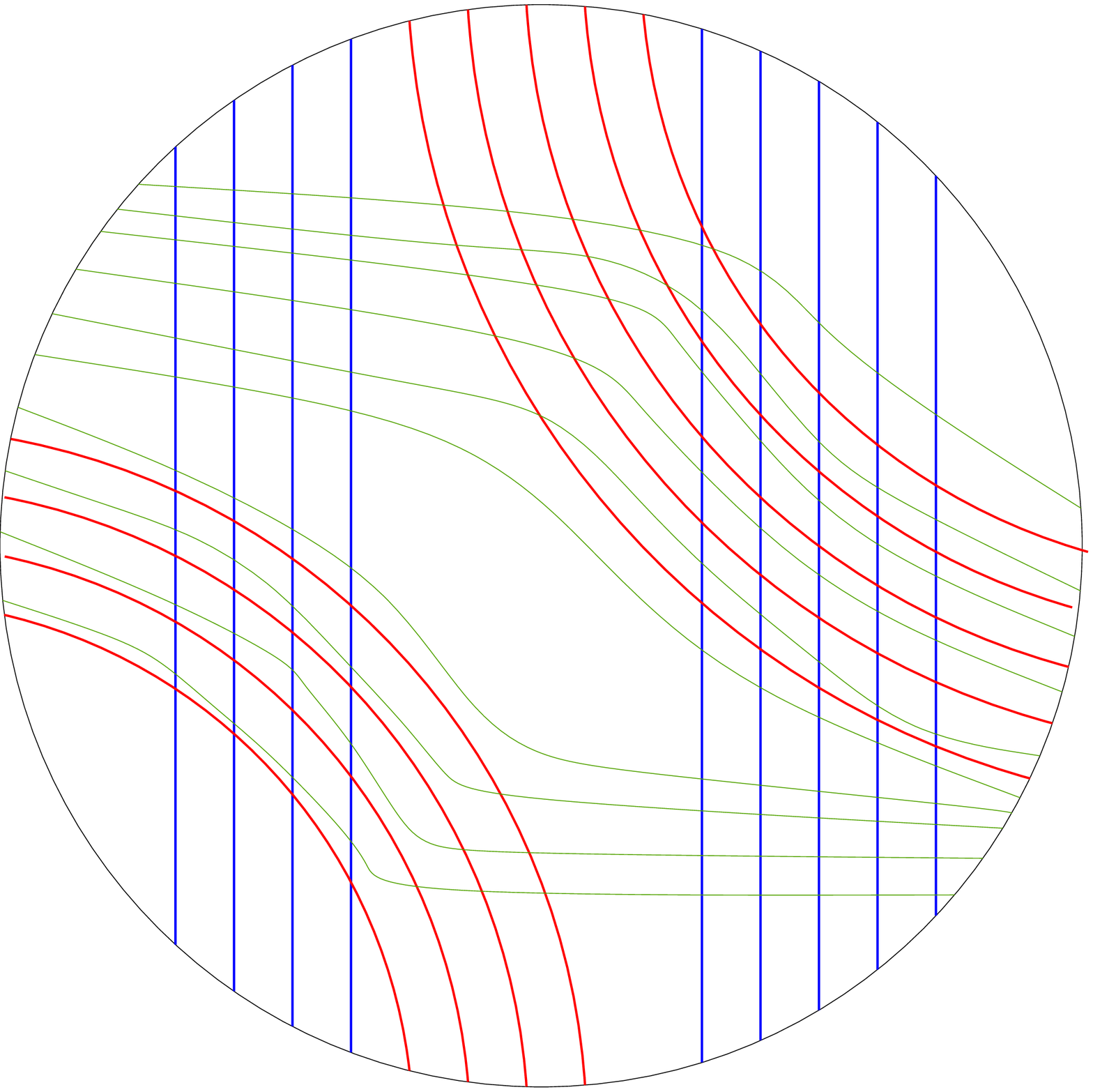}\\[0.2cm]
\mbox{front hemisphere} & \mbox{rear hemisphere}
\end{array}$
\caption{Isotoping to obtain a grid diagram}
\label{fig:isotopy}
\end{figure}
\end{proof}

\section{Proof of Theorems \ref{thm:nontransversesimple} 
}\label{sec:prooftransversenonsimple}

One way of distinguishing transverse knots in a given knot type is to prove that their $\thetah$-invariants are different. This however cannot be done straightforwardly as the vector space $\HFKh$ does not canonically correspond to a knot. 
So in order to prove that two elements are different, we have to show that there is no isomorphism of $\HFKh$ carrying one to the other. Or more explicitly it is enough to see, that there is no such isomorphism induced by a sequence of Heegaard moves. For instance if we show that one element is $0$, while the other is not, we can be certain that they are different. This is used in the proof of Theorem \ref{thm:nontransversesimple}.


\begin{proof}[Proof of Theorem \ref{thm:nontransversesimple}]
Ng, Ozsv\'ath and Thurston \cite{NgOT} showed that the knot type $10_{132}$ contains transversely non isotopic representatives $L_1$ and $L_2$ with equal self-linking number.
They proved that $\thetah(L_1)$ is zero in $\HFKh(S^3,m(10_{132}))$ while $\thetah(L_2)$ is not.
In the following we will prove, that the types $\#^n 10_{132}$ are transversely non simple. By the uniqueness of prime decomposition of knots \cite{Ad} these are indeed different knot types. Thus this list provides infinitely many examples of transversely non simple knots.
The two transversely non isotopic representatives are $\#^nL_2$ and $L_1\# (\#^{n-1}L_2)$. For the self-linking numbers we have $\textrm{sl}(\#^nL_2)=n\textrm{sl}(L_2)+(n-1)=(n-1)\textrm{sl}(L_1)+\textrm{sl}(L_2)(n-1)=\textrm{sl}(L_1\# (\#^{n-1}L_2))$.
Then we use Corollary \ref{cor:connectedsumhat} to distinguish the transverse isotopy type of $\#^nL_2$ and $L_1\# (\#^{n-1}L_2)$. There is an isomorphism from $\HFKh(S^3,m(10_{132}))\otimes\HFKh(S^3,\#^{n-1}m(10_{132}))$ to $\HFKh(S^3,\#^nm(10_{132}))$ mapping $\thetah(L_1)\otimes \thetah(\#^{n-1}L_2))=0$ to  $\thetah(L_1\#(\#^{n-1}L_2))$ thus it is zero. Similarly there is an isomorphism  mapping $\thetah(L_2)\otimes\thetah(\#^{n-1}L_2))\neq 0$ to  $\thetah(L_2\#(\#^{n-1}L_2))$ thus by induction on $n$ it does not vanish.

\end{proof}

\bibliographystyle{plain}
\bibliography{connectedsum}

\begin{thebibliography}{10}

\bibitem{BM}
J.~S. Birman and W.~M. Menasco.
\newblock Stabilization in the braid groups {II}. {T}ransversal simplicity of
  knots.
\newblock {\em Geom. Topol.}, 10:1425–--1452, 2006.
\newblock math.GT/0310280.

\bibitem{Ad}
G.~Burde and H.~Zieschang.
\newblock {\em Knots}, volume~5 of {\em de Gruyter Stud. Math.}
\newblock Walter de Gruyter, New York, 2003.

\bibitem{Ch}
Y.~Chekanov.
\newblock Differential algebra of {L}egendrian links.
\newblock {\em Invent. Math.}, 150(3):441–--483, 2002.

\bibitem{EF}
Y.~Eliashberg and M.~Fraser.
\newblock Classification of topologically trivial {L}egendrian knots.
\newblock In {\em Geometry, topology, and dynamics (Montreal, PQ, 1995)},
  volume~15 of {\em CRM Proc. Lecture Notes}, pages 17--51. Amer. Math. Soc.,
  Providence, RI, 1998.

\bibitem{EFM}
J.~Epstein, D.~Fuchs, and M.~Meyer.
\newblock Chekanov-eliashberg invariants and transverse approximations of
  {L}egendrian knots.
\newblock {\em Pacific J. Math.}, 201(1):89–--106, 2001.

\bibitem{EH}
J.~B. Etnyre and K.~Honda.
\newblock Knots and contact geometry. {I}. {T}orus knots and the figure eight
  knot.
\newblock {\em J. Symplectic Geom.}, 1(1):63–--120, 2001.
\newblock math.GT/0006112.

\bibitem{E}
J.B. Etnyre.
\newblock Legendrian and transversal knots.
\newblock In {\em Handbook of knot theory}, pages 105--185. Elsevier B. V.,
  Amsterdam, 2005.

\bibitem{EH2}
J.B. Etnyre and K.~Honda.
\newblock On connected sums and {L}egendrian knots.
\newblock {\em Adv. Math.}, 179(1):59--74, 2003.

\bibitem{MOS}
C.~Manolescu, P.~Ozsv\'ath, and S.~Sarkar.
\newblock A combinatorial description of knot {F}loer homology.
\newblock math/0607691.

\bibitem{MOSzT}
C.~Manolescu, P.~Ozsv\'ath, Z.~Szab\'o, and D.~Thurston.
\newblock On combinatorial link {F}loer homology.
\newblock math/0610559.

\bibitem{McC}
J.~McCleary.
\newblock {\em A user's guide to spectral sequences}, volume~58 of {\em
  Cambridge Studies in Advanced Mathematics}.
\newblock Cambridge University Press, Cambridge, second edition, 2001.

\bibitem{Ng}
L.~Ng.
\newblock Legendrian {T}hurston-{B}ennequin bound from {K}hovanov homology.
\newblock {\em Algebr. Geom. Topol.}, 5:1637–--1653, 2005.
\newblock math.GT/0508649.

\bibitem{NgOT}
L.~Ng, P.~Ozsv\'ath, and D.~Thurston.
\newblock Transverse knots distinguished by knot {F}loer homology.
\newblock math/0703446.

\bibitem{AAA}
P.~Ozsv\'ath and Z.~Szab\'o.
\newblock Holomorphic discs, link invariants, and the multi-variable
  {A}lexander {P}olinomial.
\newblock math/0512286.

\bibitem{OSZknot}
P.~Ozsv{\'a}th and Z.~Szab{\'o}.
\newblock Holomorphic disks and knot invariants.
\newblock {\em Adv. Math.}, 186(1):58--116, 2004.

\bibitem{OSZ3m}
P.~Ozsv{\'a}th and Z.~Szab{\'o}.
\newblock Holomorphic disks and topological invariants for closed
  three-manifolds.
\newblock {\em Ann. of Math. (2)}, 159(3):1027--1158, 2004.

\bibitem{OSzT}
P.~Ozsv\'ath, Z.~Szab\'o, and D.~Thurston.
\newblock Legendrian knots, transverse knots and combinatorial {F}loer
  homology.
\newblock math/0611841.

\bibitem{R}
J.~A. Rasmussen.
\newblock {\em Floer homology and knot complements}.
\newblock Phd thesis, Harvard University, 2003.
\newblock math/0607691.

\bibitem{SW}
S.~Sarkar and J.~Wang.
\newblock An algorithm for computing some {H}eegaard {F}loer homologies.
\newblock math/0607777.

\end{thebibliography}

\end{document}